\documentclass[graybox]{svmult}

\usepackage{type1cm}        
\usepackage{makeidx}         
\usepackage{graphicx}        
\usepackage{multicol}        
\usepackage[bottom]{footmisc}

\usepackage{newtxtext}       %
\usepackage{newtxmath}       

\usepackage{placeins}
\usepackage{float}

\usepackage{listings}
\usepackage{bbm}

\newcommand{\Nd}{N_{\delta}}

\newcommand{\bu}{\textbf{u}}
\newcommand{\bv}{\textbf{v}}
\usepackage{verbatim}
\usepackage{afterpage}
\newcommand{\og}{\omega}

\makeindex             

\graphicspath{{../images/}}

\begin{document}

\title{Weighted reduced order methods for uncertainty quantification in computational fluid dynamics}
\titlerunning{Weighted reduced order methods for UQ in CFD}
\author{Julien Genovese, Francesco Ballarin, Gianluigi Rozza and Claudio Canuto}
\institute{Julien Genovese \at Politecnico di Torino, Department of Mathematical Sciences ``G.L. Lagrange'', Corso Duca degli Abruzzi 24, 10129 Torino, Italy \and Francesco Ballarin \at Department of Mathematics and Physics, Università Cattolica del Sacro Cuore, via Garzetta 48, I-25133 Brescia, Italy \email{francesco.ballarin@unicatt.it}
	 \and Gianluigi Rozza \at MathLab, Mathematics Area, SISSA, via Bonomea 265, I-34136 Trieste, Italy \email{gianluigi.rozza@sissa.it} \and Claudio Canuto \at Politecnico di Torino, Department of Mathematical Sciences ``G.L. Lagrange'', Corso Duca degli Abruzzi 24, 10129 Torino, Italy }
%
%
\maketitle

\abstract*{In this manuscript we propose and analyze weighted reduced order methods for stochastic Stokes and Navier-Stokes problems depending on random input data (such as forcing terms, physical or geometrical coefficients, boundary conditions). We will compare weighted methods such as weighted greedy and weighted POD with non-weighted ones in case of stochastic parameters. In addition we will analyze different sampling and weighting choices to overcome the curse of dimensionality with high dimensional parameter spaces.}

\abstract{In this manuscript we propose and analyze weighted reduced order methods for stochastic Stokes and Navier-Stokes problems depending on random input data (such as forcing terms, physical or geometrical coefficients, boundary conditions). We will compare weighted methods such as weighted greedy and weighted POD with non-weighted ones in case of stochastic parameters. In addition we will analyze different sampling and weighting choices to overcome the curse of dimensionality with high dimensional parameter spaces.}

\section{Introduction}
Partial differential equations (PDEs) are often used to model complicated physical phenomena. Notable examples include the Stokes and the Navier-Stokes equations in fluid dynamics (CFD). In practical cases, the model usually depends on some parameters that can be either physical, such as the diffusivity or the maximum value of inlet boundary conditions in CFD, or geometrical, such as the length of the domain.
These parameters can be deterministic or uncertain, the latter due for example to the lack of knowledge or because of measurement errors. The model is then described by stochastic PDEs in the latter case.

As very few cases of PDEs, especially in CFD, can be solved analytically, we resort here to numerical methods to find an approximated solution. The major issue is that  classical discretization methods, to which we will refer as \textit{full order methods} in the following, take an excessive amount of computational time to e.g.\ compute moments of quantities of interest associated to stochastic PDEs. This is essentially due to expensive computational costs  related to the solution of a many-query problem, i.e.\ we need to obtain each solution associated to a (possibly large) set of parameters.

One of the methodologies devised in the last decades to bring down the computational time is to use the so called \textit{reduced order methods} (ROMs) based on a reduced basis paradigm \cite{RozzaBook,CertifiedReducedBasis}. The main idea of ROMs is to obtain a fast approximate solution for a new parameter value by using a linear combination of precomputed solutions obtained with other parameters.
In these methods we always have an \textit{offline} phase where we solve the problem several times for different parameter values, exploring the space of the parameters and storing any quantity that should be used in the next \textit{online} phase. This phase is usually computationally cheap and allows to solve the problem for another parameter without querying the underlying full order method.

When moving towards a stochastic context, some parameters are more likely to occur than others. In the setting of random variables, there is a need for fast simulations in order to perform statistics analysis on the solutions \cite{WeightedApproachPeng,rb_roxana, ch12,Stoc2,Stoc3}. 
Moreover, the probabilistic information ought to be exploited to choose the parameter values employed during the offline phase. In particular the probability can be used to (i) identify the best parameter for the offline phase of the reduced method, or (ii) as a weight for the reduced algorithms.
This work will pursue the second option, following in the footsteps of a class of \emph{weighted reduced order methods} introduced in the past \cite{WeightedApproachPeng,Stoc5,Stoc1,Stoc4} and used in many applied fields, such as optimal control and convection dominated problems \cite{CARERE2021261,torlo,zoccolan}. Two alternative offline stages will be summarized, based either on a weighted greedy algorithm \cite{WeightedApproachPeng,spann,torlo} or a weighted proper orthogonal decomposition (POD) approach \cite{LucaVenturi,Stoc2,Stoc3} . We will mainly focus on the latter, discussing several types of numerical integration: the Monte-Carlo technique, the tensor product rule and finally the sparse rule. The main contribution of this manuscript is to extend the aforementioned methods to stochastic PDEs in CFD.

The work is organized as it follows. In section 2 we introduce the strong and the weak Stokes and Navier-Stokes equations, with the related Galerkin formulation that we will use as our full order model.
In section 3 we summarize the weighted reduced order methods that we use to perform model reduction of stochastic PDEs. Section 4 will focus on tensor products and sparse grids methods, as two relevant choices to be employed in the weighted POD method. In section 5 we present some numerical results on the performance of weighted ROMs in numerical test cases in CFD. Some brief concluding remarks and future perspectives are provided in section 6.

\section{Problem setting and full order model}
We first recall some probability notions \cite{probBook} which will be used in the following.
Let us introduce a triple $(\Omega, \mathcal{F}, P)$ that denotes a complete probability space, where $\Omega$ is the set of the outcomes $\omega\in \Omega$, $\mathcal{F}$ is a $\sigma$-algebra of events and $P:\mathcal{F}\rightarrow [0,1]$ is a probability measure, i.e.\ $P(\Omega)=1$.
In such a framework we introduce a real-valued continuous \textit{random variable} $Y:(\Omega,\mathcal{F})\rightarrow (\mathbb{R},\mathcal{B}) $ being $\mathcal{B}$ the Borel $\sigma$-algebra on $\mathbb{R}$. We denote with $\Gamma$ the image of $Y$ and with $\rho:\Gamma\rightarrow \mathbb{R}$ the associated probability density function. For this random variable we can compute the \textit{expectation} of $Y$ as:
\begin{equation*}
	\mathbb{E}[Y]=\int_{\Omega}Y dP=\int_{\Gamma}y\rho(y)dy.
\end{equation*}
Higher order moments of $Y$ (e.g., the variance) can also be similarly defined.

Let us now introduce a spatial domain $D\subseteq \mathbb{R}^d$.
We define a \textit{random field} $v: D\times \Omega\rightarrow \mathbb{R}$ as a field for which $v(\mathbf{x}, \cdot)$ is a random variable for any fixed $\mathbf{x}\in D$.
We also introduce a \textit{random vector} $\textbf{v}=(v_1,v_2,...,v_d):D\times \Omega\rightarrow \mathbb{R}^d$ as a vector whose components are random fields.

We now move to the stochastic steady Stokes and Navier-Stokes problems. For this presentation we will follow \cite{Stoc1}, extending the deterministic cases of \cite{steadyStokes} and \cite{QuarteroniNPDE}.
In the following we assume a random domain $D = D(\omega): \Omega \to \mathbb{R}^d$, a random diffusivity $\nu:\Omega\rightarrow \mathbb{R}_{+}$, a random forcing term $\textbf{f}:D \times \Omega \rightarrow \mathbb{R}^d$, a random Neumann boundary condition $\textbf{h}:\partial D_N \times \Omega\rightarrow \mathbb{R}^d$ and a random inlet condition $\textbf{g}_{in}:\partial D_{in}\times \Omega \rightarrow\mathbb{R}^d$. We now introduce stochastic Stokes and Navier-Stokes problems, where the purpose is to solve the system of equations changing the outcome $\omega \in \Omega$.

For the Stokes problem we then search a pair $(\textbf{u}(\omega),p(\omega)):D \times \Omega\rightarrow \mathbb{R}^d\times \mathbb{R}$ solution to:
\begin{equation}
	\label{Ss}
	\begin{cases}
		-\nu(\omega) \Delta \textbf{u}(\textbf{x},\omega) +\nabla p(\textbf{x},\omega) =\textbf{f}(\textbf{x},\og)& \text{in }D(\og),\\
		\nabla \cdot \textbf{u}(\textbf{x},\og)=0  & \text{in $D(\og)$,}\\
		\textbf{u}(\textbf{x},\og)=\textbf{0} & \text{on $\partial D_{D,\textbf{0}}(\og)$},\\
		\textbf{u}(\textbf{x},\og)=\textbf{g}_{in}(\textbf{x},\og) & \text{on $\partial D_{in}(\og)$},\\
		\nu \dfrac{\partial \textbf{u}}{\partial \textbf{n}}(\textbf{x},\og)-p(\textbf{x},\og)\textbf{n}=\textbf{h}(\textbf{x},\og) &\text{on $\partial D_{N}(\og)$,}
	\end{cases}
\end{equation}
for almost every $\omega\in\Omega$.
Here the boundary is partitioned in $$\partial D=\partial D_{D,\textbf{0}}\cup \partial D_{in}\cup \partial D_{N}.$$

The Navier-Stokes requires the addition of a convective nonlinear term $\textbf{u}\cdot \nabla \textbf{u}$:
\begin{equation}
	\label{NS}
	\begin{cases}
		-\nu(\omega) \Delta \textbf{u}(\textbf{x},\omega) +\nabla p(\textbf{x},\omega) & \\
	\qquad \qquad  \qquad   \qquad  \quad+\textbf{u}(\textbf{x},\og)\cdot \nabla \textbf{u}(\textbf{x},\og) =\textbf{f}(\textbf{x},\og)& \text{in }D(\og),\\
		\nabla \cdot \textbf{u}(\textbf{x},\og)=0  & \text{in $D(\og)$,}\\
		\textbf{u}(\textbf{x},\og)=\textbf{0} & \text{on $\partial D_{D,\textbf{0}}(\og)$},\\
		\textbf{u}(\textbf{x},\og)=\textbf{g}_{in}(\textbf{x},\og) & \text{on $\partial D_{in}(\og)$},\\
		\nu \dfrac{\partial \textbf{u}}{\partial \textbf{n}}(\textbf{x},\og)-p(\textbf{x},\og)\textbf{n}=\textbf{h}(\textbf{x},\og) &\text{on $\partial D_{N}(\og)$.}
	\end{cases}
\end{equation}

Following \cite{WeightedApproachPeng,Stoc2}, we now move from a stochastic setting into a parametric one. In particular we suppose that $\nu,\textbf{f}$, $\textbf{h}$ and $\textbf{g}_{in}$ depend on a finite number $K$ of random variables, which we collect in a random vector $Y(\og)=(Y_1(\og),Y_2(\og),...,Y_K(\og)):\Omega \rightarrow \boldsymbol{\Gamma}=\Gamma_1\times \Gamma_2\cdots\times\Gamma_K\subset \mathbb{R}^{K}$, being $\Gamma_k$ the image of $Y_k$. We further collect each probability density function in the vector $\rho=(\rho_1,...,\rho_{K}):\Gamma\rightarrow \mathbb{R}^K$. In other words we are assuming that the stochastic dependence is expressed, with an abuse of notation, as
\begin{align*}
	& \nu(\og)=\nu(Y(\og)),\\
	& \textbf{f}(\cdot,\og)=\textbf{f}(\cdot,Y(\og)),\\
	& \textbf{h}(\cdot,\og)=\textbf{h}(\cdot,Y(\og)),\\
	& \textbf{g}_{in}(\cdot,\og)=\textbf{g}_{in}(\cdot,Y(\og)).
\end{align*}
We note that each quantity may depend only on a subset of $Y(\og)$, but we do not differentiate the stochastic dependence any further for the sake of shortness.

Thanks to the Doob-Dynkin lemma \cite{DoobDynkin} we also have that the solutions to \eqref{Ss} or \eqref{NS} satisfy:
\begin{align*}
	&\bu(\cdot, \og)=\bu(\cdot,Y(\og)),\\
	&p(\cdot,\og)=p(\cdot,Y(\og)).
\end{align*}
This result is instrumental into moving from a stochastic description towards a parametric partial differential equations. In fact, it states that we can equivalently refer to the stochastic description associated to the event $\omega\in \Omega$ or the parametric one induced by realizations of $Y(\omega)=y\in \mathbb{R}^K$. We will take the latter standpoint in the following, as it is more easily adapted to the reduced order formulation that we are willing to introduce. However, in order to do that we also need to assume that $\Gamma_k$ is a compact set for each $k$. This is indeed a very mild assumption that is often valid in numerical practice, up to a truncation of probability distributions with infinite support to a large compact set characterized by the highest probability.

The next step is to introduce a weak formulation of the problem, which requires the definition of two function spaces $V(y) = H^1_{\partial D_{D,\textbf{0}}(y)\cup \partial D_{in}(y)}(D(y); \mathbb{R}^d)$ for the velocity and $Q(y) = L^2(D(y))$ for the pressure. The bilinear and linear forms are then defined as in the deterministic case \cite{QuarteroniNPDE}:
\begin{align*}
	&a(\cdot, \cdot; y):V(y)\times V(y)\rightarrow \mathbb{R} \mbox{ such that}\\
	&a(\bu,\bv;y):=\int_{D(y)} \nu(y)\nabla \bu :\nabla  \bv \ d\textbf{x}=\sum_{i,j=1}^{d}\int_{D(y)}\nu(y) \dfrac{\partial u_i}{\partial x_j}\dfrac{\partial v_i}{\partial x_j}\ d\textbf{x},\\[20pt]
	&b(\cdot, \cdot; y):V(y)\times Q(y)\rightarrow \mathbb{R} \mbox{ such that }\\
	&b(\bv,p;y):=-\int_{D(y)}(\nabla \cdot \bv) p \ d\textbf{x}dP=-\sum_{i=1}^d \int_{D(y)}\dfrac{\partial v_i}{\partial x_i} p \ d\textbf{x},\\[20pt]
	&F(\cdot; y):V(y)\rightarrow \mathbb{R} \mbox{ such that}\\
	&F(\bv;y):=\int_{D(y)} \textbf{f}(y) \cdot \bv +\int_{\partial D_N(y)} \textbf{h}(y) \cdot \bv - a(\bu_{g_{in}}(y),\bv)=F_s(\bv;y)+F^0(\bv;y),\\[10pt]
	&\mbox{ where } F^0(\bv;y)=-a(\bu_{g_{in}}(y),\bv) \mbox{ and } F_s(\bv;y)= F(\bv;y)-F^0(\bv;y),\\[20pt]
	&G(\cdot; y):Q(y)\rightarrow \mathbb{R} \mbox{ such that }\\
	&G(q;y):=-b(\bu_{g_{in}}(y),q).
\end{align*}
Here $\bu_{g_{in}}(\cdot, y)$ is a lifting of the inlet boundary condition $\textbf{g}_{in}(\cdot, y)$. For the sake of brevity from now on we will consider $\textbf{h}=\textbf{0}$.

A cornerstone requirement to guarantee the availability of a fast reduced order model is now to trace back our domain $D(y)$ into a reference one $\hat{D}$ independent from the random vector $Y$. See e.g.\ \cite{RozzaBook} for a more detailed description. The reference domain $\hat{D}$ is often obtained as the image $D(\hat{y})$ for a specific realization $\hat{y}$. First of all we need to split our domain in a set of non-overlapping subdomains $\{D^r(y)\}_{r=1,...,R}$ such that $\cup_{r=1}^R D^r(y)=D(y)$. Next we trace back $D^r(y)$ into $\hat{D}^r$ with an affine transformation: $$\hat{\mathbf{x}}=T^r(\textbf{x}; y)=G^r(y)\textbf{x}+\textbf{g}^r(y).$$
If we introduce:
\begin{equation}
	\hat{\textbf{u}}(\hat{\mathbf{x}}; y):=\textbf{u}(T^{-1}(\hat{\mathbf{x}}; y)),\quad \hat{p}(\hat{\mathbf{x}}; y):=p(T^{-1}(\hat{\mathbf{x}}; y)),
\end{equation}
belonging to the space $\hat{V} = H^1_{\partial \hat{D}_{D,\textbf{0}}\cup \partial \hat{D}_{in}}(\hat{D}; \mathbb{R}^d)$ and $\hat{Q} = L^2(\hat{D})$, respectively,
and trace back all the bilinear and linear forms in the reference domains we obtain:
\begin{equation*}
	\hat{a}(\hat{\textbf{u}},\hat{\textbf{v}};y)=\sum_{r=1}^R \int_{\hat{D}^r} \nu(y)\sum_{i,j,j',j''}\dfrac{\partial \hat{u_i}}{\partial \hat{x}_{j'}}\big(G^r_{j'j}(y)G^r_{j''j}(y)\det(G^r(y)^{-1}) \big)\dfrac{\partial \hat{v_i}}{\partial \hat{x}_{j''}}\ d\hat{\mathbf{x}}, \mbox{ $\forall \hat{\textbf{v}} \in \hat{V}$},
\end{equation*}
and
\begin{equation*}
	\hat{b}(\hat{\textbf{v}},\hat{p};y)=-\sum_{r=1}^R\int_{\hat{D}^r}\hat{p}\big(G^r_{ij}(y)\det(G^r)^{-1}(y)\big)\dfrac{\partial \hat{v}_j}{\partial \hat{x}_i}\ d \hat{\mathbf{x}}, \mbox{ \ \ $\forall \hat{\bv}\in \hat{V},$}
\end{equation*}
and the forcing terms become:
\begin{align*}
	&\hat{F}_s(\hat{\bv};y)=\sum_{r=1}^R\int_{\hat{D}^r}\sum_{j=1}^d\big(\hat{f}^r_j \det(G^r(y))^{-1} \big)\hat{v}_j \ d \hat{\mathbf{x}},\\
	&\hat{F}^0(\hat{\bv};y) =-\hat{a}(\hat{\textbf{u}}_{\textbf{g}_{in}},\hat{\bv}),
\end{align*}
and finally:
$$
\hat{G}(\hat{q};y)=-\hat{b}(\hat{\textbf{u}}_{\textbf{g}_{in}},\hat{q}).
$$
In what follows, we will drop the hat from the bilinear forms, the unknowns and the function spaces. This slight abuse of notation is justified because from now on we will refer only to the problem traced back on the reference domain.

We can now formulate the weak form of both problems of interest, starting from the Stokes case. We want to find a solution $(\bu,p):\boldsymbol{\Gamma}\rightarrow V\times Q$ such that:
\begin{equation}
	\begin{cases}
		a(\textbf{u},\textbf{v};y)+b(\textbf{v},p;y)= F(\textbf{v};y) , \mbox{\ \ \ $\forall \textbf{v}\in V$},\\
		b(\textbf{u},q;y)=G(q;y) ,\mbox{ \ \ $\forall q \in Q$},
	\end{cases}
	\label{realstocweakstokes}
\end{equation}
for a.e. $y\in \boldsymbol{\Gamma}$ distributed according to $\rho(y)$.

In the Navier-Stokes case we only have to add the nonlinearity in its trilinear form in the first equation of \eqref{realstocweakstokes}:
$$
c(\textbf{u},\textbf{u},\textbf{v};y)=\int_{D}\sum_{i,j,m}u_i\chi_{ji}(\textbf{x};y)\dfrac{\partial u_m}{\partial x_j}v_m d\boldsymbol{x},
$$
where $\chi_{ji}$ is the $(j,i)$ element of the matrix $\chi(\textbf{x};y)=(G(\textbf{x};y)^{-1})\det(G(\textbf{x},y))$
obtaining:
\begin{equation}
	\begin{cases}
		a(\textbf{u},\textbf{v};y)+b(\textbf{v},p;y)+c(\textbf{u},\textbf{u},\textbf{v};y)= F(\textbf{v};y) , \mbox{\ \ \ $\forall \textbf{v}\in V$},\\
		b(\textbf{u},q;y)=G(q;y) ,\mbox{ \ \ $\forall q \in Q$}.
	\end{cases}
	\label{realstocNSweakstokes}
\end{equation}

We conclude this section with the full order discretization of \eqref{realstocweakstokes} and \eqref{realstocNSweakstokes} with a Galerkin approximation (e.g., the finite element method \cite{QuarteroniNPDE}). We introduce two finite dimensional spaces $V_{\Nd}$ and $Q_{\Nd}$, where, for the sake of simplicity, $\Nd$ denotes the total dimension of both $V_{\Nd}$ and $Q_{\Nd}$. So we search a pair $(\bu_{\Nd},p_{\Nd})\in V_{\Nd}\times Q_{\Nd}$ such that:
\begin{equation}
	\begin{cases}
		a(\textbf{u}_{\Nd},\textbf{v}_{\Nd};y)+b(\textbf{v}_{\Nd},p_{\Nd};y)= F(\textbf{v}_{\Nd};y),  \mbox{\ \ \ $\forall \textbf{v}_{\Nd}\in V_{\Nd}$},\\
		b(\textbf{u}_{\Nd},q_{\Nd};y)=G(q_{\Nd};y), \mbox{ \ \ $\forall q_{\Nd} \in Q_{\Nd},$}
	\end{cases}
	\label{discretestocweakstokes}
\end{equation}
for a.e. $y\in \boldsymbol{\Gamma}$, for the Stokes case, while for Navier-Stokes:
\begin{equation}
	\begin{cases}
		a(\textbf{u}_{\Nd},\textbf{v}_{\Nd};y)+b(\textbf{v}_{\Nd},p_{\Nd};y) \\ 
                   \qquad \qquad  \qquad   +c(\textbf{u}_{\Nd},\textbf{u}_{\Nd},\textbf{v}_{\Nd};y)= F(\textbf{v}_{\Nd};y),  \mbox{\ \ \ $\forall \textbf{v}_{\Nd}\in V_{\Nd}$},\\
		b(\textbf{u}_{\Nd},q_{\Nd};y)=G(q_{\Nd};y), \mbox{ \ \ $\forall q_{\Nd} \in Q_{\Nd}.$}
	\end{cases}
	\label{discretestocweakNstokes}
\end{equation}

\section{Weighted Reduced Order Methods}
We summarize now reduced order methods in a stochastic context for solving problems in computational fluid dynamics, called \textit{weighted reduced order methods}.
From the discussion in the previous section, it is now clear that the stochastic formulation and the parametrized one are formally equivalent. However, in the interest of an effective model reduction, a ROM devised on top of the parametrized formulation should properly \emph{weight} realizations of the full order system, according to the probability density function $\rho(y)$.
Several weighted reduced order methods have been devised in the past. In this paper we will follow mainly the presentation of \cite{Stoc2}, \cite{Stoc3} and \cite{Stoc1}. See for example \cite{valueatrisk} for variants with different goals, where the value at risk is used.

The goal of a ROM is to lower the computational cost from a polynomial function of $N_{\delta}$, as required by the \textit{full order methods} we introduced at the end of the previous section, to a polynomial function of $N\ll N_{\delta}$.
The methodology hinges on two different phases \cite{RozzaBook}:
\begin{itemize}
	\item \textit{offline phase}: this is the most expensive phase. Here we compute the truth solution, i.e. the solution obtained with the full order method, for several parameters chosen according to one of two algorithms called \textit{weighted greedy} and \textit{weighted POD}, we will introduce later, and according to the probability density function. The cost hence depends on a power of $\Nd$.
	In particular we want to find two reduced spaces $V_{N}$ and $Q_{N}$ with dimensions that are multiples of $N\ll\Nd$ and generated by a finite linear combination of truth solutions of the problem \eqref{discretestocweakstokes} or \eqref{discretestocweakNstokes}. For doing that we have to choose a discrete parameters space  $\mathbb{P}_h$ and searching in it some parameters $y$ for computing the truth solutions.
	
	The main novelty in the weighted approach with respect to the non-weighted one presented in \cite{RozzaBook} is that we assign different weights to the parameters according to a weight function $w(y)$, chosen following some rules related with the probability density function and discussed in the following. This function will have the role of choosing the parameters space and the basis.
	In addition in this phase we memorize several quantities in preparation for the next phase.
	
	\item \textit{online phase}: in this phase we are ready to compute fast simulations and to obtain them we search the solution in the reduced spaces using the quantities stored before. The cost of this phase must depend only on $N$ and this allows to satisfy real-time simulations and many-query problems.
\end{itemize}
What we will find in the numerical experiments is that when we are working with parameters derived from a distribution far from the uniform one (which was on the contrary the implicit assumption in the deterministic case), the weighted approach can be useful for lowering the computational cost and obtaining better approximations.

To formulate the reduced weighted method, as usual we begin assuming an \textit{affine decomposition} of the different terms in the equations. So we suppose they are a linear combination of the components of the random vector $Y(\og)=(Y_1(\og),Y_2(\og),...,Y_K(\og)):\Omega\rightarrow \mathbb{R}^K$, i.e.:
\begin{align*}
	& a(\bu,\bv;Y(\og))=a_0(\bu,\bv)+\sum_{k=1}^K a_k(\bu,\bv)Y_k(\og),\\
	& b(\bv,q;Y(\og))=b_0(\bv,q)+\sum_{k=1}^K b_k(\bv,q)Y_k(\og),\\
	& F(\bv;Y(\og))=F_0(\bv)+\sum_{k=1}^K F_k(\bv)Y_k(\og),\\
	& G(q;Y(\og))=G_0(q)+\sum_{k=1}^K G_k(q)Y_k(\og),
\end{align*}
with the bilinear forms $a_k$, $b_k$ and the linear forms $F_k$, $G_k$ are such that:
\begin{align*}
	& a_k:V\times V\rightarrow \mathbb{R}\quad \forall k=1,\dots,K,\\
	& b_k:V\times Q\rightarrow \mathbb{R}\quad \forall k=1,\dots,K,\\
	&F_k:V\rightarrow \mathbb{R}\quad \forall k=1,\dots,K,\\
	& G_k:Q\rightarrow \mathbb{R}\quad \forall k=1,\dots,K.
\end{align*}
We remind that this hypothesis, as in deterministic case, is required for a computational gain in the online phase of the reduced simulation.

Under this assumptions, the computation of pair $(\bu_{N},p_{N})\in V_{N}\times Q_{N}$ such that (i.e., the online phase):
\begin{equation}
	\begin{cases}
		a(\textbf{u}_{N},\textbf{v}_{N};y)+b(\textbf{v}_{N},p_{N};y)= F(\textbf{v}_{N};y),  \mbox{\ \ \ $\forall \textbf{v}_{N}\in V_{N}$},\\
		b(\textbf{u}_{N},q_{N};y)=G(q_{N};y), \mbox{ \ \ $\forall q_{N} \in Q_{N},$}
	\end{cases}
\end{equation}
for a.e. $y\in \boldsymbol{\Gamma}$, for the Stokes case, or
\begin{equation}
	\begin{cases}
		a(\textbf{u}_{N},\textbf{v}_{N};y)+b(\textbf{v}_{N},p_{N};y) \\ 
 \qquad \qquad  \qquad    +c(\textbf{u}_{N},\textbf{u}_{N},\textbf{v}_{N};y)= F(\textbf{v}_{N};y),  \mbox{\ \ \ $\forall \textbf{v}_{N}\in V_{N}$},\\
		b(\textbf{u}_{N},q_{N};y)=G(q_{N};y), \mbox{ \ \ $\forall q_{N} \in Q_{N},$}
	\end{cases}
\end{equation}
for the Navier-Stokes one, can be realized in a way that depends polynomially on $N$, but not on $\Nd$. See e.g.\ \cite{RozzaBook} for further details, as well as to \cite{Stokesreduced1,PODNavierStokes} for peculiarities concerning the CFD case, especially for accurate pressure recovery. We do not dwell on the online phase any more, as there are no relevant differences with respect to the deterministic case. Instead, we discuss now the two main algorithms employed during the offline phase to introduce a suitable weighting in the stochastic context.

\subsection{Weighted Greedy Algorithm}
In this part we will present the weighted greedy algorithm, following \cite{WeightedApproachPeng}, \cite{Stoc2} and \cite{Stoc3}.
First, we define a finite training set $\mathbb{P}_h\subset \boldsymbol{\Gamma}$ in which we select the snapshots to create the reduced spaces. The training set is chosen according to the probability density
function $\rho(y)$.
At the beginning we choose a parameter $y_1 \in \mathbb{P}_h$ using the probability distribution $\rho(y)$, we solve the truth problem \eqref{discretestocweakstokes} or \eqref{discretestocweakNstokes} and we obtain the solution $\textbf{U}_{\Nd}(y_1)=(\bu_{\Nd}(y_1),p_{\Nd}(y_1))$. We note that with this approach we need to treat the solution as a pair. We define the reduced spaces $V_N:=\text{span}(\bu_{\Nd}(y_1))$ and $Q_N:=\text{span}(p_{\Nd}(y_1))$, with $N=1.$ After this initialization the algorithm proceeds iteratively: at the $n$-th iteration we search the parameter $y_{n}$ such that $\textbf{U}_{\Nd}(y_n):=(\textbf{u}_{\Nd}(y_n),p_{\Nd}(y_n))$ is the worst approximated solution by $\textbf{U}_{N}(y_n):=(\textbf{u}_{N}(y_n),p_{N}(y_n))$ in $\mathbb{P}_h$. Then, we compute the truth solution with this parameter and we add it to the reduced space, i.e., $V_n:=\text{span}(\lbrace \bu_{\Nd}(y_k)\rbrace_{k=1}^n)$ and $Q_n:=\text{span}(\lbrace p_{\Nd}(y_k)\rbrace_{k=1}^n)$.

In the deterministic setting, the greedy approach relies on the choice of  $y\in \mathbb{P}_h$ such that:
\begin{equation}
	\mbox{arg}\max_{y\in \mathbb{P}_h} \Delta_{N}^{\Nd}(y),
\end{equation}
where  $\Delta_{N}^{\Nd}(y)$ is the error estimator (see e.g.\ \cite{RozzaBook} for its expression) and, defined $W:=V\times Q$, it is such that
\begin{equation}
	||\textbf{U}_{\Nd}(y)-\textbf{U}_{N}(y)||_W \leq \Delta_{N}^{\Nd}(y).
\end{equation}
The rationale is to select the parameter that maximizes the error between the truth solution and the reduced one. A weighted variant can thus be introduced to take into account the probabilistic information on the parameters. Let $w:\boldsymbol{\Gamma}\rightarrow \mathbb{R}$ be a weight function (to be defined later), and introduce the weighted error estimator $\hat{\Delta}_{N}^{\Nd}(y)$ as
\begin{equation}
	\hat{\Delta}_{N}^{\Nd}(y):= w(y) \ \Delta_{N}^{\Nd}(y).
\end{equation}
The weighted greedy algorithm thus seeks
\begin{equation}
	\mbox{arg}\max_{y\in \mathbb{P}_h} \hat{\Delta}_{N}^{\Nd}(y).
\end{equation}
We note that following this approach we weigh two different contributions in the definition of the iterative greedy procedure: (i) parameters $y$ that result in a large error (estimator) are favored, and (ii) parameters which are more likely for the probability density function are favored.

The greedy algorithm drives the weighted error estimator $\hat{\Delta}_{N}^{\Nd}(y)$ (i.e.\, an upper bound to the weighted error $||\textbf{U}_{\Nd}(y)-\textbf{U}_{N}(y)||_Y w(y)$) to decrease. Among the possible choices for the weight $w$, we mention:
\begin{itemize}
	\item $w(y):=\sqrt{\rho(y)}$ can be set if we desire the greedy algorithm to control the mean of $||\textbf{U}_{\Nd}-\textbf{U}_{N}||_W^2$. Indeed
	\begin{align*}
		&\mathbb{E}\bigg[||\textbf{U}_{\Nd}-\textbf{U}_{N}||^2_W\bigg]=\int_{\Omega}||\textbf{U}_{\Nd}-\textbf{U}_{N}||_W^2 dP =\\
		& \int_{\boldsymbol{\Gamma}} ||\textbf{U}_{\Nd}(y)-\textbf{U}_{N}(y)||_W^2\rho(y)\ dy \leq \int_{\boldsymbol{\Gamma}}\Delta_{N}^{\Nd}(y)^2\rho(y)\ dy,
	\end{align*}
	so that
	\begin{equation}
		\mathbb{E}[||\textbf{U}_{\Nd}-\textbf{U}_{N}||_W^2]\leq \int_{\boldsymbol{\Gamma}}\hat{\Delta}_{N}^{\Nd}(y)^2 \ dy\leq |\boldsymbol{\Gamma}|\sup_{y\in \boldsymbol{\Gamma}}\hat{\Delta}_{N}^{\Nd}(y)^2,
	\end{equation}
	with $|\boldsymbol{\Gamma}|< +\infty$ being the finite measure of the set $\boldsymbol{\Gamma}$, note that we have supposed the compactness of $\boldsymbol{\Gamma}$.
	\item $w(y):=\rho(y)$, which leads to the control
	\begin{equation}
		||\mathbb{E}[\textbf{U}_{\Nd}]-\mathbb{E}[\textbf{U}_{N}]||_W\leq \int_{\boldsymbol{\Gamma}}||\textbf{U}_{\Nd}(y)-\textbf{U}_{N}(y)||_W\rho(y)\ dy\leq |\boldsymbol{\Gamma}| \sup_{y\in \boldsymbol{\Gamma} }\hat{\Delta}_{N}^{\Nd}(y).
	\end{equation}
\end{itemize}

We finish this section presenting the weighted greedy algorithm:\\
Initialization:
\begin{itemize}
	\item take a discrete space $\mathbb{P}_h$ according to the probability density function
	\item take a tolerance $\epsilon_{tol}$ as stopping criteria for the algorithm
	\item choose a maximum number of reduced bases $N_{max}$
	\item choose a first parameter $y_1$ and create the sample space $S_1:=\{y_1\}$
	\item solve the truth problem for $y_1$ and create $W_{1}:=\mbox{span}\{\textbf{U}_{\Nd}(y_1)\}$
\end{itemize}
	Iteration:
\begin{itemize}
	\item for $N=2,\dots,N_{max}$
	\begin{itemize}
		\item choose $y_N\in \mathbb{P}_h$ such that it maximizes $\hat{\Delta}_{N-1}^{\Nd}(y)$
		\item if $\hat{\Delta}_{N-1}^{\Nd}(y_N)\leq \epsilon_{tol}$, then $N_{max}=N$, else
		\item solve the truth problem for $y_N$ to obtain $\textbf{U}_{\Nd}(y_N)$
		\item add $y_N$ to the space $S_{N-1}$, creating $S_{N}=S_{N-1}\cup \{y_N\}$
		\item add  $\textbf{U}_{\Nd}(y_N)$ to the reduced basis space $W_{N-1}$
		\item creating  $W_{N}:=W_{N-1}\oplus\mbox{span}\{\textbf{U}_{\Nd}(y_N)\}$
	\end{itemize}
\end{itemize}
We finally note that with this approach we need to compute only $N$ truth solutions for the reduced space and this is an advantage with respect to the POD approach that we will see in the next section. In contrast in this case we have to know the error estimator that it is problem-dependent and sometimes it is not available. The weighted POD approach which we summarize in the following has actually been proposed to overcome this limitation.
\subsection{Weighted Proper Orthogonal Decomposition Algorithm (POD)}
In the \textit{weighted POD approach} we would like to find the $N$-dimensional subspaces $V_{N}$ and $Q_N$ such that the first one minimizes the error:
\begin{equation}
	\int_{\boldsymbol{\Gamma}}||\bu_{\Nd}(y)-\bu_{N}(y)||^2_V\rho(y)\ dy,
	\label{stochasticintegral}
\end{equation}
while the second one:
\begin{equation}
	\int_{\boldsymbol{\Gamma}}||p_{\Nd}(y)-p_{N}(y)||_{Q}^2\rho(y)\ dy.
\end{equation}
where $\textbf{u}_N$ and $p_N$ are the projections of the \textit{truth solutions} $\bu_{\Nd}$ and $p_{\Nd}$ over the reduced spaces, i.e. $\textbf{u}_N=P_{\textbf{u}}(\textbf{u}_{N_{\delta}})$ and $p_N=P_p(p_{N_{\delta}})$ where $P_{\textbf{u}}$ and $P_{p}$ are the projection operators over $V_{N}$ and $Q_N$, respectively. As in the weighted greedy approach, a weight (equal to the probability density function) has been embedded in the definition of the energy in order to incorporate the stochastic dependence.

Since the treatment is analogous for the two spaces, from now on we will refer only to the first one.
Let us start discretizing the integral \eqref{stochasticintegral} with a finite sum, so the goal is now to minimize:
\begin{equation}
	\sum_{y\in \mathbb{P}_h}w(y)||\bu_{\Nd}(y)-\bu_{N}(y)||_V^2=\sum_{i=1}^M w_i||\bu_{\Nd}(y_i)-P_{N}(\bu_{\Nd}(y_i))||_V^2,
	\label{finitestocdiscretization}
\end{equation}
where $\mathbb{P}_h$ is a finite discretization of $\boldsymbol{\Gamma}$ of cardinality equal to $M$.
The points $y_i$ in $\mathbb{P}_h$, as well as the weights $w_i$, have to be chosen according to some quadrature rule. In the next section we will see examples such as the Monte-Carlo method, the tensor product rule or the Smolyak rule.

In order to minimize the quantity \eqref{stochasticintegral} we can follow a similar strategy  to the one used in the deterministic case \cite{Stokesreduced1} and \cite{Stokesreduced2}. We introduce the operator $\hat{C}_{\textbf{u}}$ defined as:
$$
\hat{C}_{\textbf{u}}(\bv_{\Nd}):=\sum_{m=1}^M w_m(\bv_{\Nd},\psi_m)_V\psi_m,
$$
where $\psi_m:=\bu_{\Nd}(y_m)$ and $\bv_{\Nd}\in V_M:=\mbox{span}\{\bu_{\Nd}(y_i)| y_i\in  \mathbb{P}_h\}$. We note that for creating this operator we need $M$ truth solutions.
Afterwards, in order to minimize the quantity \eqref{finitestocdiscretization}, we search the leading eigenfunctions and eigenvalues of $\hat{C}_{\textbf{u}}$. Following the same derivation as in the deterministic case \cite{RozzaBook} we arrive to the algebraic formulation in which we have to find the eigenvectors of $\hat{C}_{\textbf{u}}:=P\cdot C_{\textbf{u}}$ where $C_{\textbf{u}}$ is the matrix defined as:
$$
C_{\textbf{u}}:=S_{\textbf{u}}^TX_{\textbf{u}}S_{\textbf{u}}\in \mathbb{R}^{M\times M},
$$
where, as in \cite{Stokesreduced1} and \cite{Stokesreduced2}:
$$
S_{\textbf{u}}:=\bigg[\textbf{u}(y_1)|\textbf{u}(y_2)| ...|\textbf{u}(y_M)\bigg]\in \mathbb{R}^{N_\textbf{u}\times M}\quad (X_\textbf{u})_{ij}:=\bigg(\boldsymbol{\phi}_j,\boldsymbol{\phi}_i \bigg)_V,
$$
where $\{\boldsymbol{ \phi}_i\}_{i=1}^{\Nd}$ are the Lagrangian basis of the truth problem \cite{QuarteroniNPDE}. We proceed similarly for the pressure.

Finally $P$ is defined as $$P:=\mbox{diag}(w_1,...,w_M),$$ and it plays the role of a preconditioner matrix.

We note that $\hat{C}_{\textbf{u}}$ is not symmetric in the usual sense but it is with respect to the scalar product induced by the matrix $C$. In fact, if we consider the induced scalar product defined as:
$$
\langle x,y\rangle_C:=x^{T}Cy,
$$
$\hat{C}$ is symmetric with respect to this scalar product if and only if
\begin{equation}
	\langle \hat{C}x,y \rangle_C= \langle x, \hat{C}y \rangle_C,
\end{equation}
that, using the definition, is equivalent to
$$
x^T\hat{C}^TCy=x^TC\hat{C}y,
$$
and so $\hat{C}$ is symmetric if
$$
\hat{C}^TC=C\hat{C}.
$$
Indeed,
$$
\hat{C}^TC=(PC)^TC=C^TP^TC=CPC=C\hat{C},
$$
for the symmetry of $C$ and $P$.
So with this scalar product we can use the spectral theorem and we have an orthogonal basis of eigenvectors.

Let us summarize now the \textit{weighted POD algorithm}:
\begin{itemize}
	\item choose the training set $\mathbb{P}_h\subset \boldsymbol{\Gamma}$ and the weights $w_i$ according to some
	quadrature rule
	\item solve the truth problem for each of the parameters
	in $\mathbb{P}_h$ and find the solutions $\{\psi_i\}_{i=1,\dots,M}$
	\item assemble the matrix $\hat{C}_{ij}=w_iC_{ij}$ with $C_{ij}=(S_{\textbf{u}}    ^TX_{\textbf{u}}S_{\textbf{u}})_{ij}$
	and search the $N$ biggest eigenvectors $\xi_1,\dots,\xi_{N}$ and related eigenvalues
	\item construct $V_{N}=\mbox{span}\{\xi_1,\dots,\xi_{N}\}$.
\end{itemize}

\section{Tensor products and sparse grids methods}
In this section we will treat in more details the choice of the set $\mathbb{P}_h$ and the associated weights resulting from a quadrature rule. We will follow the work of \cite{Smolyakarticle,SmolyakQuadrature}.

The topic of the quadrature rules is usually seen in a more general context, as approximation of a multidimensional integral of a function $f: \mathbb{R}^d \rightarrow \mathbb{R}$ with a weight function $g:\mathbb{R}^d\rightarrow \mathbb{R}$ as follows:
\begin{equation*}
	\begin{split}
	\mathcal{I}_g^d f=& \int_{\boldsymbol{\Gamma}} g(y)f(y)dy = \\
	& \int_{\Gamma_1\times \Gamma_2 \times ...\times \Gamma_d}g(y^{(1)},y^{(2)},...,y^{(d)})f(y^{(1)},y^{(2)},...,y^{(d)})\,\,dy^{(d)} ... dy^{(1)},
	\end{split}
\end{equation*}
where
$\Gamma_i\subset \mathbb{R}$ is an interval and so $\boldsymbol{\Gamma}:=\Gamma_1\times \Gamma_2 \times \dots \times \Gamma_d$ is an hyper-rectangle in $\mathbb{R}^d$, $$g(y^{(1)},y^{(2)},...,y^{(d)})=g_1(y^{(1)})\cdot \cdot \cdot g_d(y^{(d)})$$ is the weight function where $g_i:\Gamma_i\rightarrow \mathbb{R}$, $g_i(y^{(i)})\geq 0$ $\forall i$: in particular in our case it is the probability density function of a random vector $Y$ with probabilistic independent components.
There are several methods for approximating this integral, we focus on the following:
\begin{itemize}
	\item Monte-Carlo method;
	\item tensor product rule \cite{numericalintegration};
	\item sparse Smolyak rule \cite{Smolyakarticle}.
\end{itemize}
The Monte-Carlo is the easiest way to do this integration but it has a slow convergence rate to the truth integral \cite{Smolyakarticle,SmolyakQuadrature} and together with the tensor product rule they suffer of the problem of the curse of dimensionality while the sparse Smolyak rule tries to overcome it.

With Monte-Carlo algorithm we interpret the equation \eqref{stochasticintegral} as a mean in the probability space $\Omega$:
\begin{equation}
	\mathbb{E}\big[||\bu_{\Nd}(Y)-\bu_{N}(Y)||^2_V\big].
\end{equation}
So we can use the theory used for a general function $f:\Omega\rightarrow\mathbb{R}$
\begin{equation*}
	\mathbb{E}\big[f(Y)\big] ,
\end{equation*}
where $Y$ is random variable with density distribution $\rho$ with compact support.

So in general we have to compute:
\begin{equation*}
	\mathbb{E}\big[f(Y)\big] = \int_{J_1\times J_2 \times ...\times J_d}\rho(y^{(1)},y^{(2)},...,y^{(d)})f(y^{(1)},y^{(2)},...,y^{(d)})\,\,dy^{(d)} ... dy^{(1)}.
\end{equation*}
The resulting Monte-Carlo method approximates such expectation in the following way:
\begin{equation}
	\mathbb{E}\big[f(Y)\big]\approx\sum_{i=1}^{N}\dfrac{1}{N}f(Y_i),
\end{equation}
where $Y_i=y_i\in \boldsymbol{\Gamma}$ is a realization of the random variable. So in our case we can choose $\mathbb{P}_h$ taking $M$ realization of the random variable $Y$ and $w(y)=\dfrac{1}{M}$.

The two remaining approaches instead select $\mathbb{P}_h$ and $w$ according to some quadrature rule for approximating \eqref{stochasticintegral}. So if we take a quadrature rule $\mathcal{Q}_{\rho}$, considering the function $\rho$ in the integral and an integrable function $f:\boldsymbol{\Gamma}\rightarrow \mathbb{R}$ we have:
\begin{equation}
	\mathcal{Q}_{\rho}(f):=\sum_{i=1}^{M}w_i f(y_i),
\end{equation}
where $\{y_i\}$ and $\{w_i\}$ are a set of nodes and associated weights that approximates the general integral:
\begin{equation}
	\int_{\Omega}f(y)\rho(y)\ dy.
\end{equation}
If we change $\mathcal{Q}_{\rho}$ obviously we change the nodes, the weights and so the preconditioner $P$.

In tensor product quadrature rules we take first a set of univariate quadrature rules $(U_k^{(j)})_{j=1}^d$ where $k$ is the number of nodes used for the approximation.
In our case the rule is chosen depending on $g_j$ according to the numerical integration \cite{numericalmethods}.
So for each $j=1,\dots,d$ we have a set of $k$ nodes $(y_i^{(j)})_{i=1}^k$ and weights $(w_i^{(j)})_{i=1}^k$ associated with the rule $U_k^{(j)}$.
So we can approximate the integral in the following way:
\begin{equation}
	\mathcal{I}^d_g f \approx \sum_{i_1=1}^{k}\cdot \cdot \cdot \sum_{i_d=1}^k w_{i_1}^{(1)}\cdot \cdot \cdot w_{i_d}^{(d)} f(y_{i_1}^{(1)},...,y_{i_d}^{(d)})=
	\bigotimes_{i=1}^dU^{(i)}_k(f).
	\label{tensorProduct}
\end{equation}
as a tensor product quadrature of order $k$.

\begin{figure}[ht]
	\centering
	\includegraphics[scale=0.4]{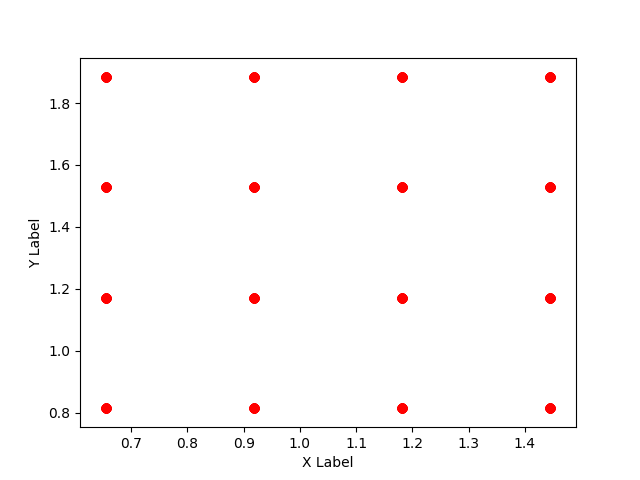}
	\caption{Tensor product grid with $4$ nodes for a \textit{Gauss-Jacobi} univariate rule ($g_j(y_j)=(1-y_j)^{\alpha}(1+y_j)^{\beta}$) with $d=2$.d}
	\label{TensorProductGrid}
\end{figure}

We can see an example of tensor product grid using a Gauss-Jacobi univariate quadrature rule in figure \ref{TensorProductGrid}.

With this grid we have a problem with \textit{the curse of dimensionality}. In fact if we use $n$ nodes for each dimension, at the end we have $n^d$ nodes. So if we increase the dimension of the space of the parameters, the number of nodes involved will grow exponentially.

\begin{figure}[ht]
	\centering
	\advance\leftskip-3cm
	\advance\rightskip-3cm
	\includegraphics[scale=0.4]{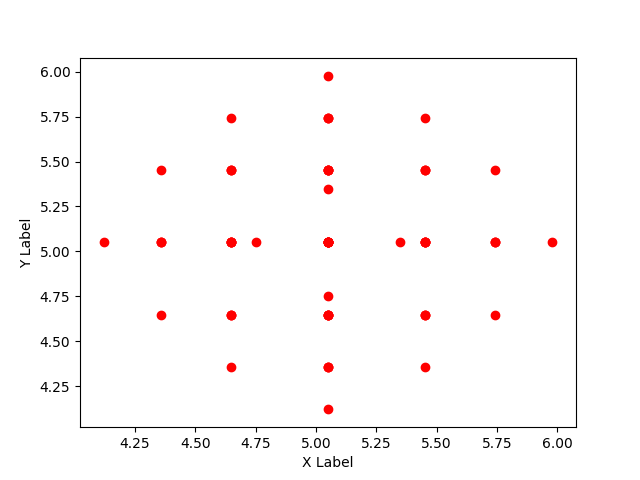}
	\caption{Sparse tensor grid $\mathcal{Q}^2_6$ with a \textit{Gauss Jacobi} rule.}
	\label{Q26}
\end{figure}

Now let us pass to the sparse grid to overcome this issue.
\newtheorem{mydef}{Definition}
\begin{mydef}
	(Smolyak quadrature rule): let $(U_i^{(j)})_{i=1}^\infty$ be a sequence of univariate quadrature rules in the interval $\emptyset \not =J_j \subset \mathbb{R}, j=1,...,d$, where $U^{(j)}_i$ has $i$ nodes and weights.
	
	We introduce the \textit{difference operators} in $J_j$ by setting
	\begin{align*}
		\Delta_0^{(j)} = 0,  && \Delta_1^{(j)}= U_1^{(j)},  &&\Delta_{i+1}^{(j)}=U_{i+1}^{(j)}-U_{i}^{(j)}, \mbox{  for $i\in \mathbb{N}_0$}.
	\end{align*}
	The \textit{Smolyak quadrature rule} of order $k$ in the hyper-rectangle $J_1\times J_2\times \cdot \cdot \cdot  \times J_d$ is the operator
	\begin{equation}
		\mathcal{Q}_k^d=\sum_{|\alpha|_1\leq k , \alpha\in \mathbb{N}^d} \bigotimes_{j=1}^{d}\Delta_{\alpha_j}^{(j)}.
	\end{equation}
\end{mydef}
We now recall a theorem proved in \cite{miatesi}:
\begin{theorem}
	Let $\alpha\in \mathbb{N}^d$ and $\alpha\geq \mathbbm{1}$, i.e. $\alpha_j\geq 1, \forall j$. Then
	\begin{equation}
		\bigotimes_{j=1}^{d}\Delta_{\alpha_j}^{(j)}=\sum_{\substack{\gamma \in \{0,1\}^d
				\alpha-\gamma\geq \mathbbm{1}}}(-1)^{|\gamma|_1}\bigotimes_{j=1}^{d} U_{\alpha_j-\gamma_j}^{(j)}.
	\end{equation}
\end{theorem}
This theorem is important to write the quadrature rule in a new way:
\begin{equation}
	\mathcal{Q}_k^d=\sum_{|\alpha|_1\leq k , \alpha\in \mathbb{N}^d} \bigotimes_{j=1}^{d}\Delta_{\alpha_j}^{(j)}=\sum_{|\alpha|_1\leq k , \alpha\in \mathbb{N}^d} \sum_{\substack{\gamma \in \{0,1\}^d
			\alpha-\gamma\geq \mathbbm{1}}}(-1)^{|\gamma|_1}\bigotimes_{j=1}^{d} U_{\alpha_j-\gamma_j}^{(j)}.
	\label{sparse}
\end{equation}
We can see an example of the grid associated with the quadrature rule in figure \ref{Q26}.

What is important to note from \eqref{sparse} is that the sparse rule does a combination of tensor product quadrature with lower order with respect to those in \eqref{tensorProduct} trying to obtain the same order of approximation but requiring less nodes.
\FloatBarrier

\section{Numerical Experiments}
In this section we will present some numerical experiments to assess the numerical reliability of the proposed weighted ROMs. In particular we will treat the stochastic Stokes and Navier-Stokes problems with
\begin{align*}
	& \nu(y)=1,\\
	& \textbf{f}(\cdot,y)=0,\\
	& \textbf{h}(\cdot,y)=0,\\
	& \textbf{g}_{in}(\mathbf{x},y) = (- v_{max}\, x_1 \cdot (x_1 - 3), 0).
\end{align*}

\begin{figure}[ht]
	\centering
	\advance\leftskip-3cm
	\advance\rightskip-3cm
	\includegraphics[scale=0.3]{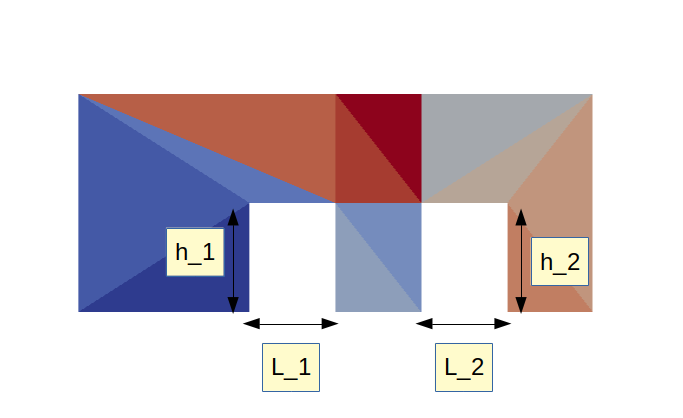}
	\caption{Geometry of the problem with four geometrical dependencies: $L_1,h_1,L_2,h_2$.}
	\label{figure}
\end{figure}

The geometry is that one in figure \ref{figure}, where the triangles denotes the decomposition of our domain ($\Omega=\cup_{r=1}^R \Omega^r$) as explained before.
For what concerns the boundaries, the inlet section $\partial D_{in}$ is on the left side, the outlet $\partial D_{N}$ is on the right side, while the remaining boudaries form a rigid wall.

The parameters involved are:
\begin{equation}
	y=(y^{(0)},y^{(1)},y^{(2)},y^{(3)},y^{(4)})=(L_1,h_1,L_2,h_2,v_{max}).
\end{equation}
The first four parameters are geometrical and referred to in figure \ref{figure}. With $v_{max}$ we denote the maximum value of the parabola in the boundary condition

The ranges considered for the five parameters are:
\begin{align*}
	&y^{(0)}\in (0.2, 1.9),\\
	&y^{(1)}\in(0.2, 2.0),\\
	&y^{(2)}\in(0.2, 1.9),\\
	&y^{(3)}\in(0.2, 2.0),\\
	&y^{(4)}\in(0.2, 20.0).
\end{align*}
We have taken around $M=240$ parameters for the training set $\mathbb{P}_h$ in all the experiments. Sometimes, with the tensor product rule and the Smolyak rule, we have obtained more or less parameters depending on the case because it is more difficult to have the exact number that we want due to the complexity of the construction.

In all cases the parameters are taken from a Beta distribution with density distribution (rescaled to the required range):
\begin{equation}
	f(y;\alpha,\beta):=\dfrac{1}{B(\alpha,\beta)}y^{\alpha-1}(1-y)^{\beta-1},
\end{equation}
with $\dfrac{1}{B(\alpha,\beta)}=\dfrac{\Gamma(\alpha+\beta)}{\Gamma(\alpha)\Gamma(\beta)}$ where $\Gamma$ is the Gamma function.

\begin{figure}[ht]
	\centering
	\advance\leftskip-3cm
	\advance\rightskip-3cm
	\graphicspath{ {../images/} }
	\includegraphics[scale=0.35]{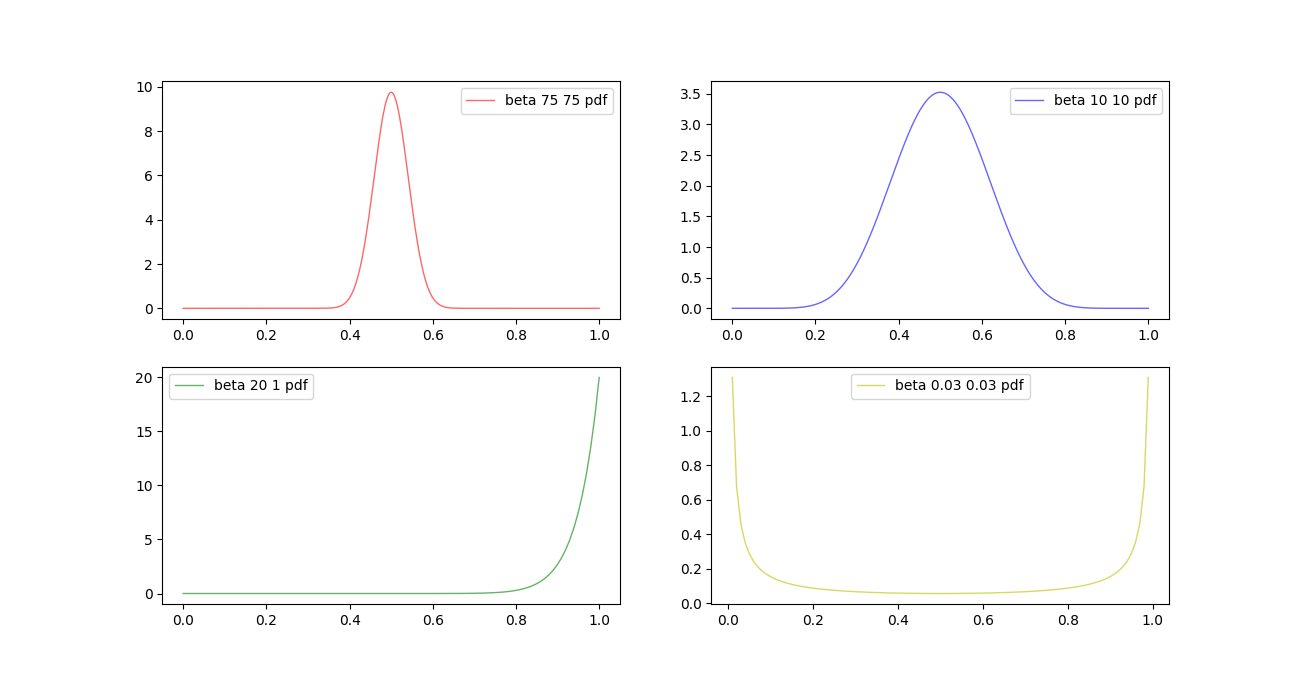}
	\caption{Some types of Beta distributions changing $\alpha$ and $\beta$. Top row: from left to right, $(\alpha, \beta) = (75, 75)$ and $(\alpha, \beta) = (10, 10)$. Bottom row: from left to right, $(\alpha, \beta) = (20, 1)$ and $(\alpha, \beta) = (0.03, 0.03)$.}
	\label{beta}
\end{figure}

Changing the values of $\alpha$ and $\beta$ we can change completely the shape of the distribution as we can see in figure \ref{beta} where we have reported some examples.
The numerical experiments have been done changing each time the probability distribution and the algorithm used.

For the Stokes case we have tried in a case to compare several stochastic and deterministic approaches: a standard greedy approach, a weighted greedy one, a standard POD one and finally a weighted POD one. In another case we have compared some numerical integration methods for the POD case: a Monte-Carlo POD, a tensor product POD and a sparse rule POD.

In the Navier-Stokes case we have tried to compare a standard POD approach with a weighted POD one. We also tried a Monte-Carlo POD, a tensor product POD and a sparse rule POD as in the Stokes case to see if we obtain the same conclusions.

In the plots that we will propose we will see the absolute error and the relative error with a $H^1-$semi-norm for the velocity on the $y-$axis while on the $x-$axis we have the number $N$ of basis used for the reduced solution. They are plotted with a logarithmic scale. In addition we will show the maximum error and the relative maximum error with the $L^{\infty}$-norm.

These errors have been computed taking $100$ parameters obtained randomly according to the chosen distribution. For each one we have found the truth solution and the reduced one changing the number of basis.

Mathematically, for each $N$ we have computed:
\begin{align*}
	&\mbox{ absolute error: } \int_{\Omega}|\bu_{\Nd}(Y(\og))-\bu_{N}(Y(\og))|_{H^1}dP\approx \frac{1}{100}\sum_{i=1}^{100} |\bu_{\Nd}(y_i)-\bu_{N}(y_i)|_{H^1},\\
	&\mbox{ absolute maximum error: } \max_{i=1,...,100}|\bu_{\Nd}(y_i)-\bu_{N}(y_i)|_{H^1},\\
	& \mbox{ relative error: } \frac{1}{100}\dfrac{\sum_{i=1}^{100} |\bu_{\Nd}(y_i)-\bu_{N}(y_i)|_{H^1}}{|\bu_{\Nd}(y_i)|_{H^1}},\\
	&\mbox{ relative maximum error: } \max_{i=1,...,100}\dfrac{|\bu_{\Nd}(y_i)-\bu_{N}(y_i)|_{H^1}}{|\bu_{\Nd}(y_i)|_{H^1}}.
\end{align*}
All the following computations have been done using RBniCS library \cite{RBniCS} which is based on FEniCS \cite{logg2012automated}.

In the Stokes case we will focus in particular on two types of Beta distributions, one with $(\alpha,\beta)=(0.03,0.03)$ and one with $(\alpha,\beta)=(75,75)$. As we can see in figure \ref{beta} they are rather different because the first one is concentrated in two zones while the second one in only one. In \cite{miatesi} we can find other experiments with  $(\alpha,\beta)=(10,10)$, so with the same shape of $(\alpha,\beta)=(75,75)$ but less concentrated, and $(\alpha,\beta)=(20,1)$ always concentrated in one zone but not symmetric.
\begin{figure}[ht]
	\centering
	\advance\leftskip-3cm
	\advance\rightskip-3cm
	\includegraphics[scale=0.4]{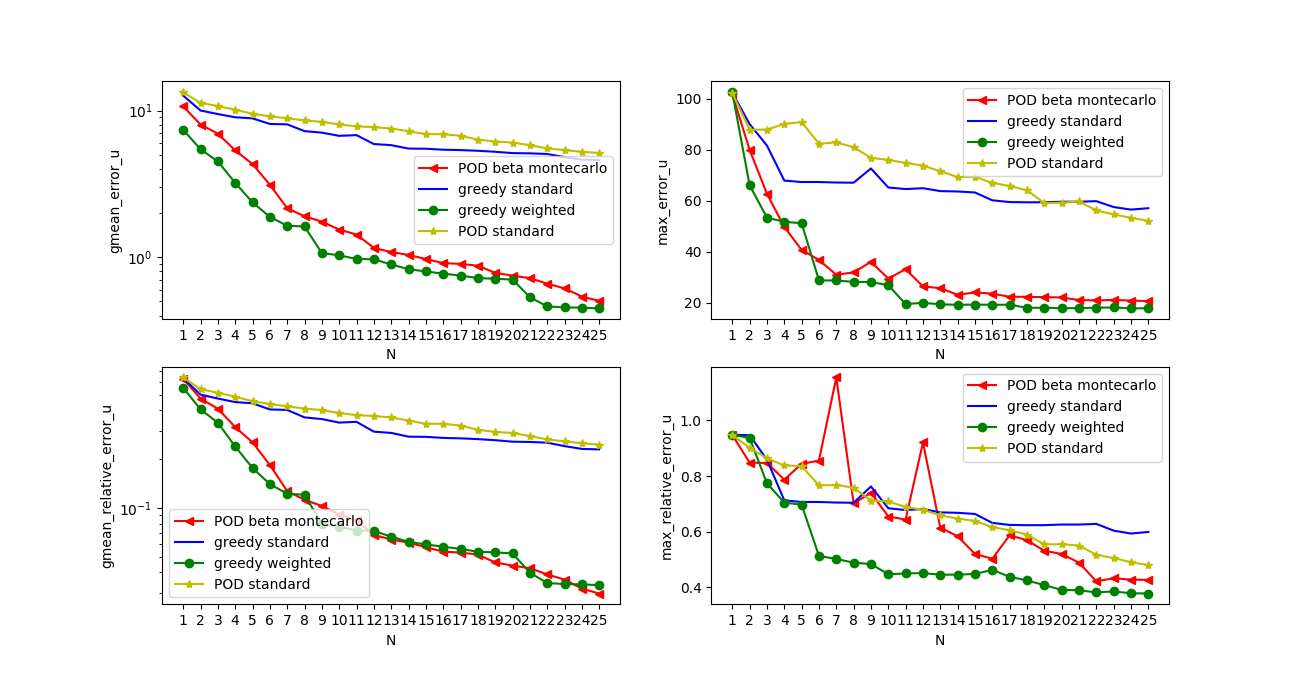}
	\caption{Stokes: standard greedy, weighted greedy, standard \textit{POD} and weighted \textit{POD}, with a $Beta(0.03,0.03).$}
	\label{s_velocity1_003}
\end{figure}
\begin{figure}[ht]
	\centering
	\advance\leftskip-3cm
	\advance\rightskip-3cm
	\includegraphics[scale=0.4]{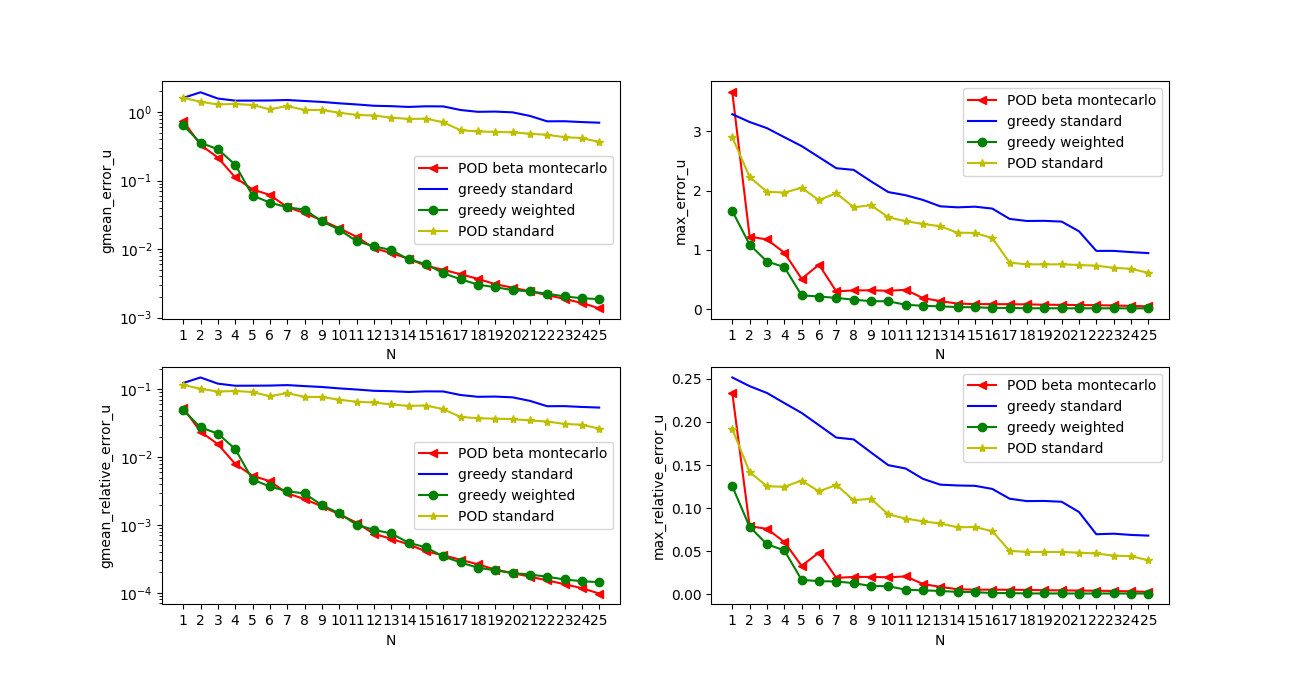}
	\caption{Stokes: standard greedy, weighted greedy, standard \textit{POD} and weighted \textit{POD}, with a $Beta(75,75)$.}
	\label{s_velocity1_75}
\end{figure}
\begin{figure}[ht]
	\centering
	\advance\leftskip-3cm
	\advance\rightskip-3cm
	\includegraphics[scale=0.4]{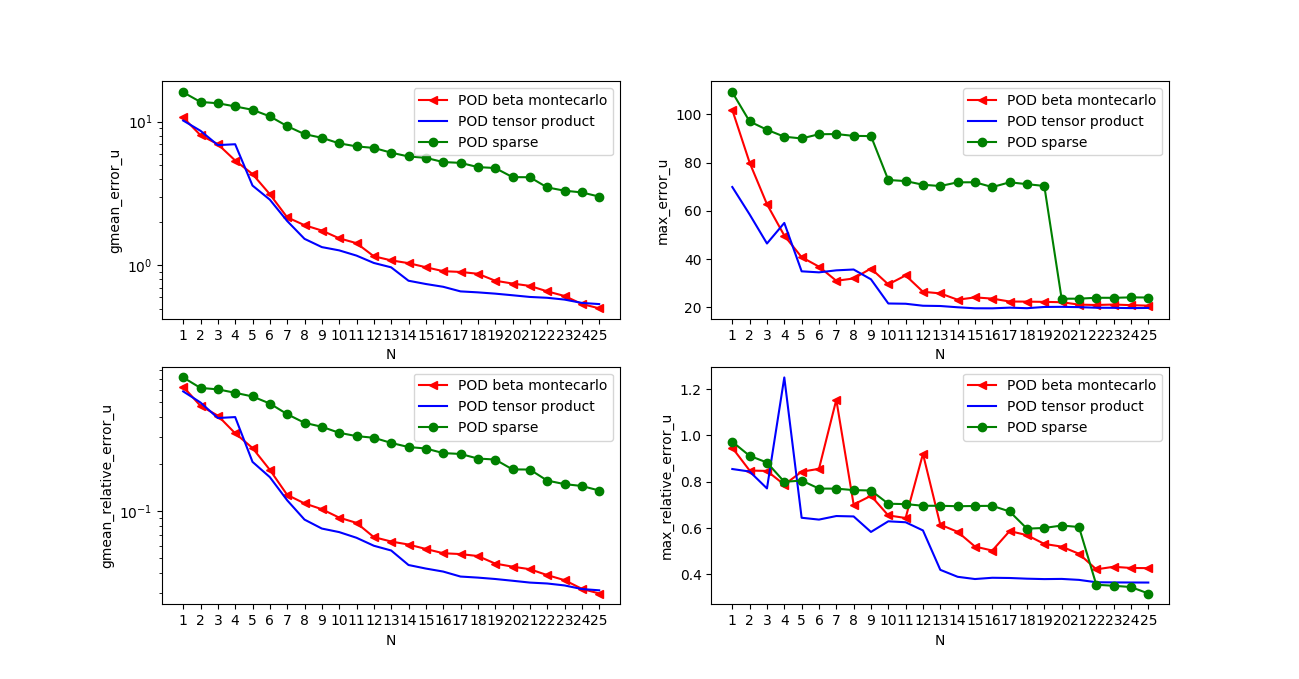}
	\caption{Stokes: Monte-Carlo \textit{POD}, Tensor product \textit{POD} with a Gauss-Jacobi quadrature rule, Sparse rule \textit{POD}, with a $Beta(0.03,0.03)$.}
	\label{s_velocity2_003}
\end{figure}
\begin{figure}[ht]
	\centering
	\advance\leftskip-3cm
	\advance\rightskip-3cm
	\includegraphics[scale=0.4]{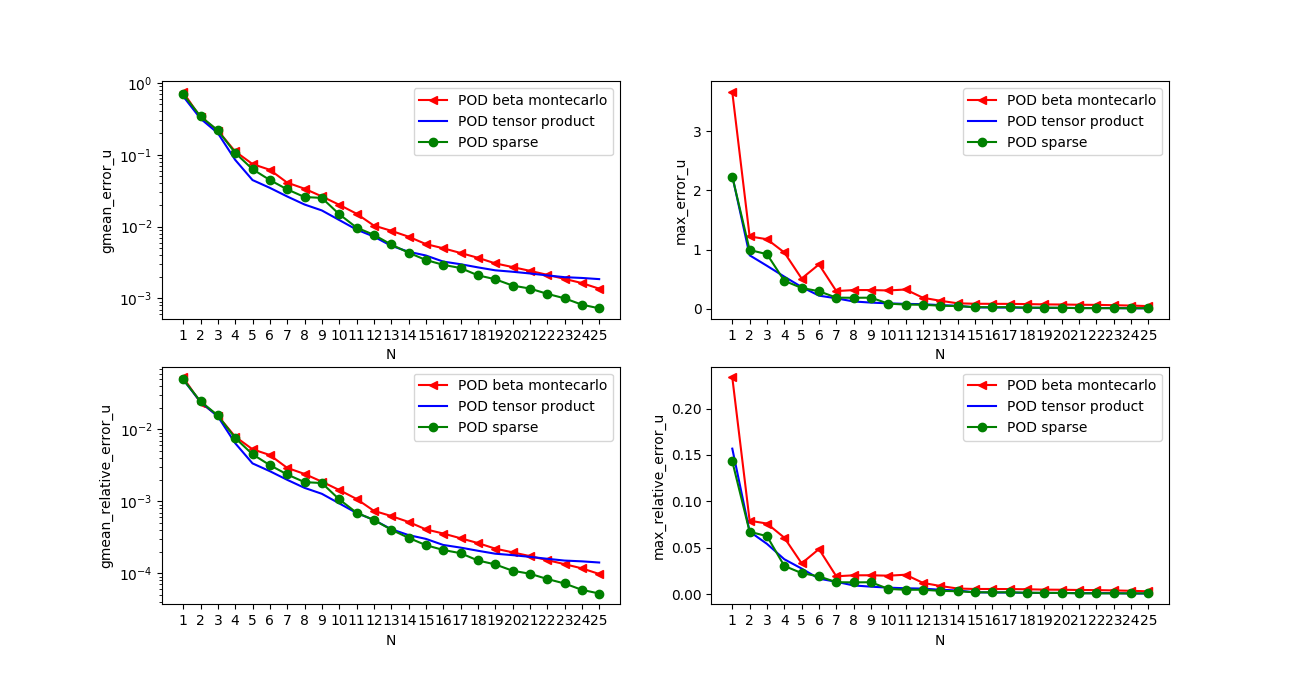}
	\caption{Stokes: Monte-Carlo \textit{POD}, Tensor product \textit{POD} with a Gauss-Jacobi quadrature rule, Sparse rule \textit{POD}, with a $Beta(75,75)$.}
	\label{s_velocity2_75}
\end{figure}

Let us begin to analyse the Stokes problem. As we can see from figures \ref{s_velocity1_003} and \ref{s_velocity1_75} the weighted methods are better that the non-weighted ones but it seems that there are no big differences between a weighted greedy approach and a weighted POD one. We can also note that the accuracy is much less in the case of $Beta(0.03,0.03)$ with respect to $Beta(75,75)$. This is probably due to the fact that in the first case the probability density function is concentrated in two zones and not in only one as in the $Beta(75,75)$ and so it is more difficult to obtain a good accuracy with few basis functions. We see the same problem in figures \ref{s_velocity2_003} and \ref{s_velocity2_75}. In fact we see that the sparse rule using a more complex combination of parameters is less accurate with respect to the other methods in $Beta(0.03,0.03)$. In the case of $Beta(75,75)$ is likely easier to chose the parameters and so all the methods have the same accuracy.

\begin{figure}[ht]
	\centering
	\advance\leftskip-3cm
	\advance\rightskip-3cm
	\includegraphics[scale=0.4]{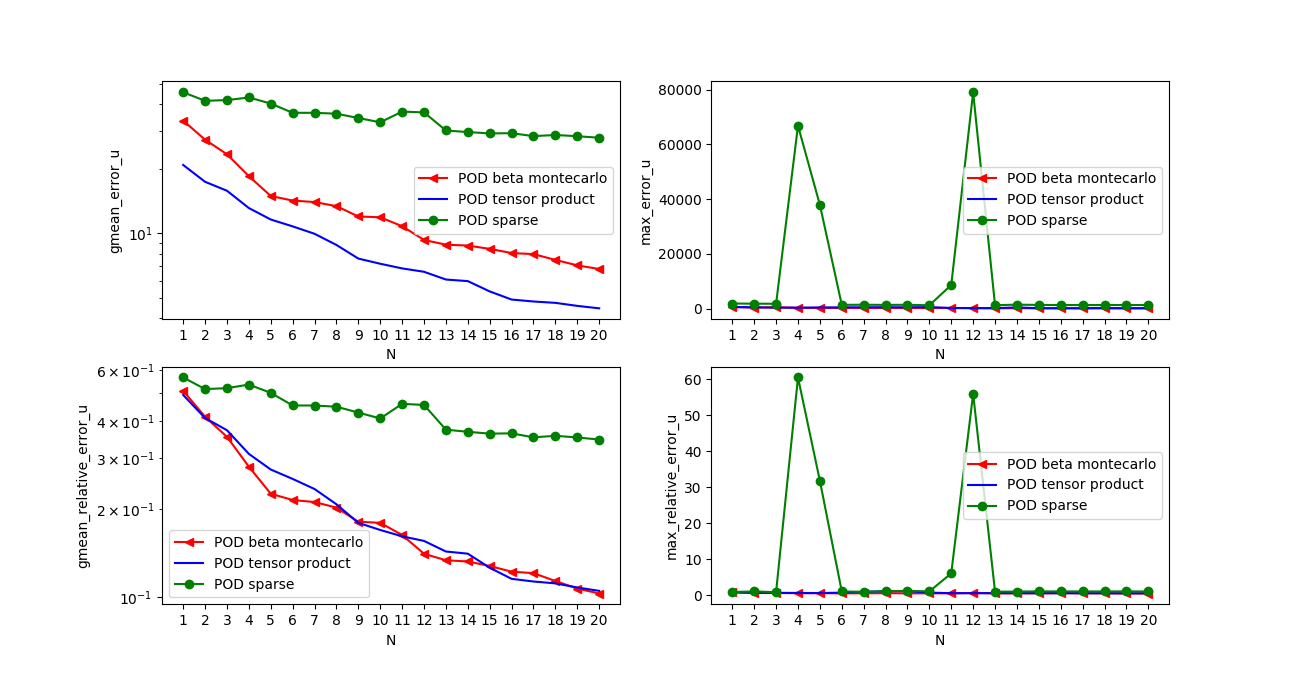}
	\caption{Navier-Stokes: Monte-Carlo \textit{POD} and tensor product rule \textit{POD} with a Gauss-Jacobi quadrature rule, sparse rule \textit{POD}, with a $Beta(0.03,0.03)$.}
	\label{ns_velocity2_003}
\end{figure}

\begin{figure}[ht]
	\centering
	\advance\leftskip-3cm
	\advance\rightskip-3cm
	\includegraphics[scale=0.4]{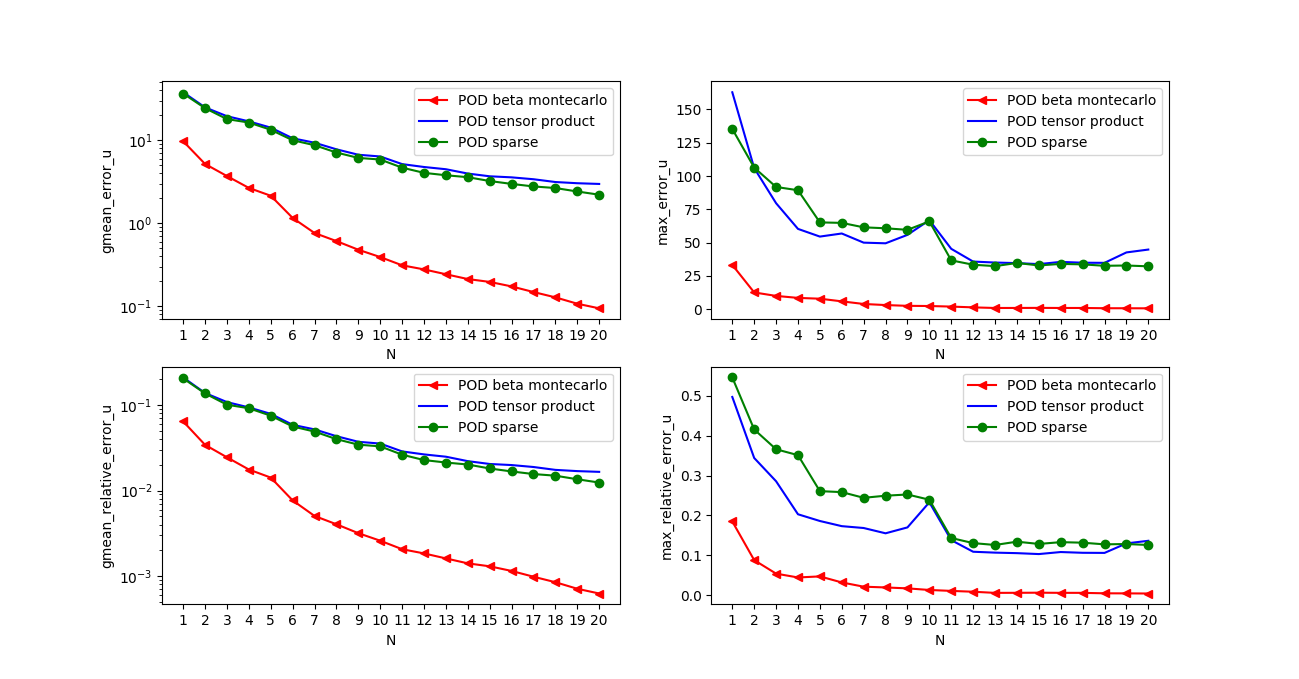}
	\caption{Navier-Stokes: Monte-Carlo \textit{POD} and tensor product rule \textit{POD} with a Gauss-Jacobi quadrature rule, sparse rule \textit{POD}, with a $Beta(10,10)$.}
	\label{ns_velocity2_10}
\end{figure}

\begin{figure}[ht]
	\centering
	\advance\leftskip-3cm
	\advance\rightskip-3cm
	\includegraphics[scale=0.4]{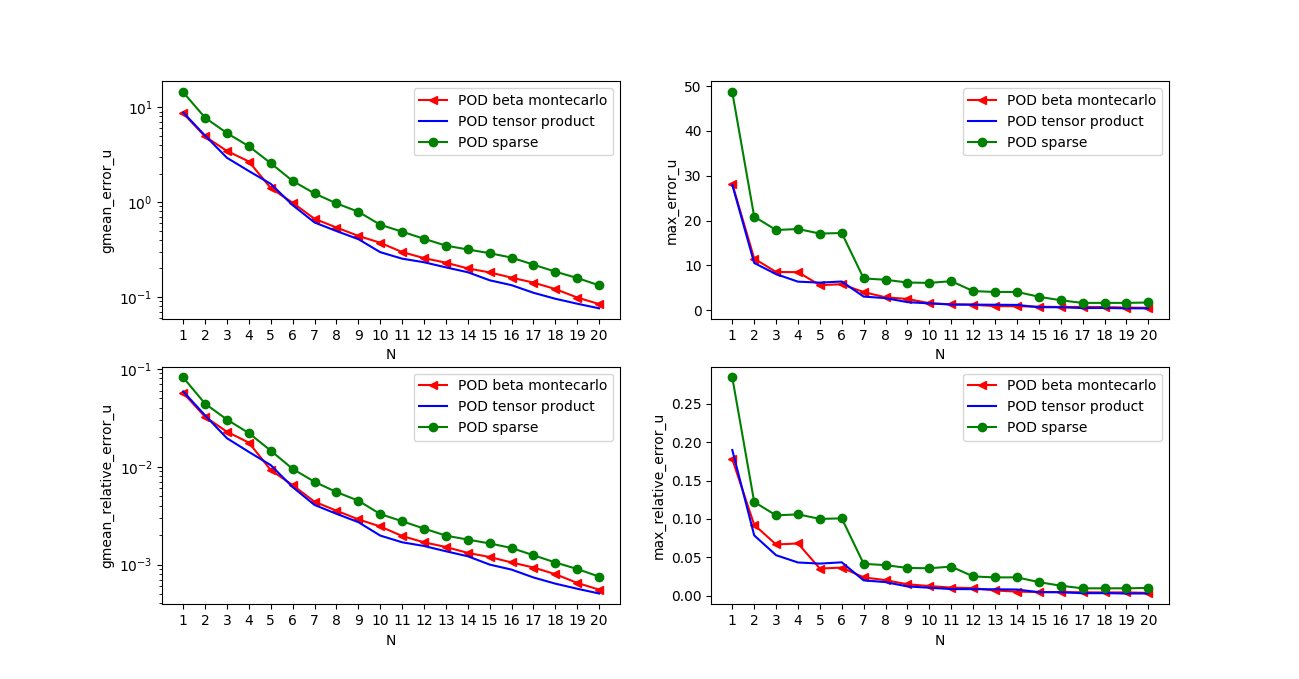}
	\caption{Navier-Stokes: Monte-Carlo \textit{POD} and tensor product rule \textit{POD} with a Gauss-Jacobi quadrature rule, sparse rule \textit{POD}, with a $Beta(75,75)$.}
	\label{ns_velocity2_75}
\end{figure}

\FloatBarrier
Let us pass now to the Navier-Stokes case where we also introduce an experiment with a $Beta(10,10)$.
In this case we have proposed in the first three figures \ref{ns_velocity2_003}, \ref{ns_velocity2_10} and \ref{ns_velocity2_75} the comparison between the three weighted POD methods.
We can observe that as in the Stokes case is not easy to obtain a good accuracy with a $Beta(0.03,0.03)$ in particular for the sparse rule. Since the training sample of the sparse rule is chosen away from the most probable area, the accuracy on the test set is quite low. The other techniques show more reliable results. For the $Beta(10,10)$ the POD with Monte Carlo sampling outperforms the other quadrature rules. Then, when we work with a $Beta=(75,75)$ all the methods obtain good and similar results.

\begin{figure}[ht]
	\centering
	\advance\leftskip-3cm
	\advance\rightskip-3cm
	\includegraphics[scale=0.4]{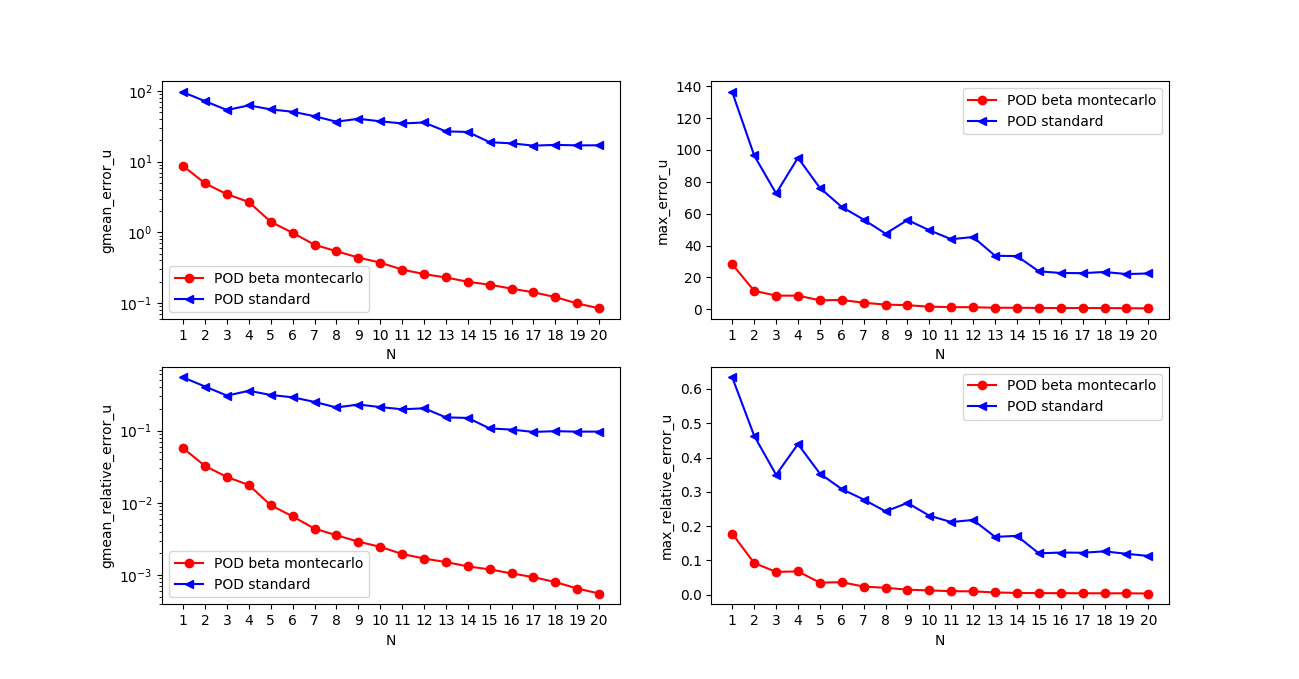}
	\caption{Navier-Stokes, first experiment: standard \textit{POD} and weighted Monte-Carlo \textit{POD}, with a $Beta(75,75)$.}
	\label{ns_velocity1_75}
\end{figure}

Finally, analyzing figure \ref{ns_velocity1_75} we can see that also in this case the weighted method heavily outperforms the non-weighted one. The deterministic method also has a very slow decay of the error with respect to $N$.

In these numerical experiments we have obtained similar results for the Stokes case and the Navier-Stokes one. As expected, we have seen that the weighted algorithms work better than the standard ones, in particular for the case of concentrated probability distributions, e.g. the $Beta(75,75)$.

We also have verified that the sparse Smolyak rule is reliable and so it is a good choice to avoid the curse of dimensionality and to obtain good results for the numerical integration, but in the cases where the distribution is concentrated in more than one zone we need more parameters that in a one-zone concentrated distribution.

The tensor product rule does not give particular advantages with respect to the Monte-Carlo one but it is more complicated to implement, so we think it is not a good solution for our problem.
In general, Monte Carlo sampling according to the underlying probability is the best solution with respect to all the other possibilities when we do not have too many parameters.

\FloatBarrier

\section{Conclusions and perspectives}
In this paper we have proposed weighted reduced basis methods for Stokes and Navier-Stokes problems with stochastic terms. We have first introduced the strong and the weak formulations of the two problems, followed by the reduced ones. In particular we have treated the weighted greedy algorithm and the weighted POD one. This last one has some variants according to Monte-Carlo rule, tensor product rule and sparse rule.
In our experiments we have tried to compare the different POD algorithms trying to see if the sparse rule can be an alternative against the curse of dimensionality.
We have seen that it is a good rule when the probability is concentrated in a small interval close to the center of the parameter domain. Nevertheless, when we treat problems with a probability distribution that it is not concentrated in one zone, we need more reduced basis to obtain a good approximation of the truth solution.

For possible future investigations we think that the study of other nonlinear stochastic problems can be interesting, as for the case of a nonlinear elastic beam \cite{artcaso1} or the nonlinear Schr\"odinger equation \cite{artcaso2}, especially in case of large parametric dimensions. Weighted variants of hyperreduction techniques \cite{artcaso3} need to be employed in those cases.
We finally think that other types of sparse grids can be implemented, as described in \cite{IUQ} or with a more general approach in \cite{SmolyakQuadrature}.
This might be beneficial especially in industrial problems characterized by several uncertain parameters.

\section*{Acknowledgements}
We acknowledge the support by European Union Funding for Research and Innovation -- Horizon 2020 Program -- in the framework of European Research Council Executive Agency: Consolidator Grant H2020 ERC CoG 2015 AROMA-CFD project 681447 ``Advanced Reduced Order Methods with Applications in Computational Fluid Dynamics''. We also acknowledge the MIUR PRIN 2017  ``Numerical Analysis for Full and Reduced Order Methods for the efficient and accurate solution of complex systems governed by Partial Differential Equations'' (NA-FROM-PDEs) and the INDAM-GNCS project ``Tecniche Numeriche Avanzate per Applicazioni Industriali''.
The computations in this work have been performed with RBniCS \cite{RBniCS} library, developed at SISSA mathLab, which is an implementation in FEniCS of several reduced order modelling techniques; we acknowledge developers and contributors to both libraries.

\FloatBarrier

\bibliographystyle{abbrv}
\bibliography{references.bib}

\end{document}